\pdfoutput=1

\documentclass{article}

\usepackage{graphicx}
\usepackage{amsfonts,amssymb}
\usepackage{amsmath, amsthm}

\usepackage{fullpage}

\newtheorem{theorem}{Theorem}
\newtheorem{lemma}{Lemma}

\theoremstyle{remark}
\newtheorem{remark}{Remark}

\def\diver{{\rm div}}
\def \nn {\mathbf{n}}
\def \RR {\mathbb R}
\def\cT{{\mathcal T}}
\def\nn{\mathbf{n}}
\def\dx{\,\d x }
\def\dS{\,\d S }

\def\grad{\nabla\!}
\def\bn{\mathbf{n}}

\def\LdK{{L^2(\partial K)}}
\def\LK{{L^2(K)}}

\def\pd{\partial}
\def\d{\mathrm{d}}

\def\qed{$\hfill\square$}
\def \Hppj{{H^{p+1}}}

\def \Th {{\mathcal T}_h}

\def \sumK {\sum_{K\in\cT_h}}
\def \intK {\int_K}

\def\intdK{\int_{\partial K}}
\def \<{\langle}
\def \>{\rangle}

\def \vt{\vartheta}

\def \tuh{\tilde{u}_h}
\def \teh{\tilde{e}_h}

\def \pdO{\partial\Omega}
\def \tu{\tilde{u}}
\def \intdKa{\int_{\partial K^-\setminus\partial\Omega}}
\def \intdKb{\int_{\partial K^-\cap\partial\Omega}}
\def \emx{e^{\mu(x,t)}}
\def \emmx{e^{-\mu(x,t)}}
\def \edmx{e^{2\mu(x,t)}}
\def \emdmx{e^{-2\mu(x,t)}}
\def \emdmxK{e^{-2\mu(x_K,t)}}

\def \txi{\tilde{\xi}}
\def \intO{\int_\Omega}
\def \vt{\vartheta}
\def \tPi{{\Pi}}
\def \tx{{\tilde{x}}}
\def \tt{{\tilde{t}}}

\begin{document}

\title{On the time growth of the error of the DG method for advective problems}

\author{%
{\sc
V\' aclav Ku\v cera\thanks{Corresponding author. Email: vaclav.kucera@email.cz}} \\[2pt]
Charles University, Faculty of Mathematics and Physics\\
Sokolovsk\'{a} 83, Praha 8, 186\,75, Czech Republic\\[6pt]
{\sc and}\\[6pt]
{\sc Chi-Wang Shu}\thanks{Email: shu@dam.brown.edu}\\[2pt]
Division of Applied Mathematics, Brown University\\
Providence, RI 02912, USA.
}

\maketitle

\begin{abstract}
{In this paper we derive a priori $L^\infty(L^2)$ and $L^2(L^2)$ error estimates for a linear advection-reaction equation with inlet and outlet boundary conditions. The goal is to derive error estimates for the discontinuous Galerkin (DG) method that do not blow up exponentially with respect to time, unlike the usual case when Gronwall's inequality is used. While this is possible in special cases, such as divergence-free advection fields, we take a more general approach using exponential scaling of the exact and discrete solutions. Here we use a special scaling function, which corresponds to time taken along individual pathlines of the flow. For advection fields, where the time massless particles carried by the flow spend inside the spatial domain is uniformly bounded from above by some $\widehat{T}$, we derive $O(h^{p+1/2})$ error estimates where the constant factor depends only on $\widehat{T}$, but not on the final time $T$. This can be interpreted as applying Gronwall's inequality in the error analysis along individual pathlines (Lagrangian setting), instead of physical time (Eulerian setting).
} 
\end{abstract}

\section{Introduction}
The discontinuous Galerkin (DG) finite element method introduced in \cite{rh} is an increasingly popular method for the numerical solution of partial differential equations. The DG method was first formulated for a neutron transport equation and such problems remain the major focus of the DG community. Since it is for problems of advective or convective nature that
the DG method is suited best and shows its strengths as opposed to other numerical methods for such problems.

In this paper, we shall consider a scalar time-dependent linear advection-reaction equation of the form 
\begin{equation}
\frac{\pd u}{\pd t}+a\cdotp\grad\, u+cu=0.
\label{intro:contprob}
\end{equation}
We will discretize the problem in space using the discontinuous Galerkin method with the upwind numerical flux on unstructured simplicial meshes in $\RR^d$. The first analysis of the DG method for a stationary advection-reaction problem with constant coefficients was made in \cite{Lesaint1974}, later improved in \cite{Johnson}. 

In this paper we derive a priori estimates of the error $e_h=u-u_h$, where $u_h$ is the DG solution. The goal is to derive estimates of $e_h$ in the $L^\infty(L^2)$ and $L^2(L^2)$ norms of the order $Ch^{p+1/2}$, where the constant $C$ does not depend exponentially on the final time $T\to +\infty$. Such results already exist in the literature, however they are derived under the ellipticity condition
\begin{equation}
c-\tfrac{1}{2}\diver a\ge\gamma_0>0
\label{intro:elipt_cond}
\end{equation}
for some constant $\gamma_0$. In the paper \cite{Svadlenka} $\gamma_0=0$ is also admissible, which corresponds to the interesting case of divergence-free advection field $a$. Without assuming (\ref{intro:elipt_cond}), when one proceeds straightforwardly, at some point Gronwall's inequality must be used in the proofs. This however results in exponential growth of the constant factor $C$ in the error estimate with respect to time.

Here, we will circumvent condition (\ref{intro:elipt_cond}) by considering an exponential scaling transformation 
\begin{equation}
u(x,t)=e^{\mu(x,t)}\tu(x,t),\quad u_h(x,t)=e^{\mu(x,t)}\tu_h(x,t)
\label{intro:expon_scaling}
\end{equation}
of the exact and discrete solutions. Substituting (\ref{intro:expon_scaling}) into (\ref{intro:contprob}) and its discretization results in a new equation for $\tu$ with a modified reaction term $\tilde{c}=\frac{\pd \mu}{\pd t}+a\cdotp\grad\mu+c$. The function $\mu$ can then be chosen in various ways in order to satisfy the ellipticity condition (\ref{intro:elipt_cond}) for the new equation. After performing the error analysis, the resulting error estimates for $\teh=\tu-\tuh$ can be transformed back via (\ref{intro:expon_scaling}) to obtain estimates for $e_h$.

The use of transformations similar to (\ref{intro:expon_scaling}) is not new. The trick known as \emph{exponential scaling}, corresponding to taking $\mu(x,t)=\alpha t$ for some sufficiently large constant $\alpha$, is very well known from the partial differential equation and numerical analysis communities. However, the application of such a scaling transformation inherently leads to exponential dependence of the error on $t$ after transforming from $\teh$ back to $e_h$. Other choices of $\mu$ are possible, e.g. $\mu(x,t)=\mu_0\cdotp x$ for some constant vector $\mu_0$. This choice was used in the stationary case in \cite{navert} and \cite{Johnson}. Taking $\mu(x)$ independent of $t$ corresponds to the analysis of the stationary case from \cite{Marini}. In the nonstationary case the use of this stationary transform would lead to too restrictive assumptions on the advection field.

In this paper we construct the function $\mu$ from (\ref{intro:expon_scaling}) using characteristics of the advection field. Namely, $\mu$ will be proportional to time along individual pathlines of the flow. If pathlines exist only for a finite time bounded by some $\widehat{T}$ for each particle entering and leaving the spatial domain $\Omega$, we obtain error estimates exponential in $\widehat{T}$ and not the final physical time $T$. This is a result we would expect if we applied Gronwall's lemma in the Lagrangian framework along individual pathlines (which exist only for a finite time uniformly bounded by $\widehat{T}$) instead of the usual application of Gronwall with respect to physical time in the Eulerian framework. The analysis can be carried out under mild assumptions that $a(\cdot,\cdot)$ satisfies the assumptions of the Picard-Lindel\"{o}f theorem and that there are no characteristic boundary points on the inlet (i.e. $a\cdot\bn$ is uniformly bounded away from zero on the inlet).

The paper is organized as follows. After introducing the continuous problem in Section \ref{sec:contprob}, we discuss the exponential scaling transform, its variants and application  in the weak formulation in Section \ref{sec:exp_scaling}. In Section \ref{sec:DG} we introduce the DG formulation and its basic properties. In Sections \ref{sec:Anal_adv_react_terms} and \ref{sec:MethLines} we analyze the DG advection and reaction forms and estimate the error of the method. Finally, in Section \ref{sec:mu_construct} we construct the scaling function $\mu$ and prove necessary results on its regularity and other properties needed in the analysis.

We use $(\cdot,\cdot)$ to denote the $L^2(\Omega)$ scalar product and $\|\cdot\|$ for the $L^2(\Omega)$ norm. To simplify the notation, we shall drop the argument $\Omega$ in Sobolev norms, e.g. $\|\cdot\|_{H^{p+1}}$ denotes the $H^{p+1}(\Omega)-$norm. We will also denote the Bochner norms over the whole considered interval $(0,T)$ in concise form, e.g. $\|u\|_{L^2(H^{p+1})}$ denotes the $L^2(0,T;H^{p+1}(\Omega))$-norm of $u$. Throughout the paper, $C$ will be a generic constant independent of $h$.

\section{Continuous problem}
\label{sec:contprob}
Let $\Omega\subset\RR^d$, $d\in\mathbb{N}$ be a bounded polygonal (polyhedral) domain with Lipschitz boundary $\pd\Omega$. Let $0<T\leq+\infty$ and let $Q_T=\Omega\times(0,T)$ be the space-time domain. We note that we admit the time interval to be infinite in the special case of $T=+\infty$.

We seek $u:Q_T\to\RR$ such that
\begin{alignat}{2}
\frac{\pd u}{\pd t}+a\cdotp\grad\, u+cu&=0 &&\text{in }Q_T,\label{contprob1}\\
u&=u_D &&\text{on }\pd\Omega^-\times(0,T),\label{contprob2}\\
u(x,0)&=u^0(x),\quad  && x\in\Omega\label{contprob3}.
\end{alignat}
Here $a:\overline{Q_T}\to\RR^d$ and $c:\overline{Q_T}\to\RR$ are the given advective field and reaction coefficient, respectively. By $\pdO^-$ we denote the \emph{inflow} part of the boundary, i.e.
\begin{equation}
\pdO^-=\{x\in\pdO; a(x,t)\cdot\nn(x)<0, \forall t\in(0,T)\},
\nonumber
\end{equation}
where $\nn(x)$ is the unit outer normal to $\pdO$ at $x$. For simplicity we assume that the inflow boundary remains the same, but in general $\pd\Omega^-$ can depend on $t$. We do not consider this case since it would only make the notation more complicated without changing the analysis itself. 

We assume that the reaction coefficient satisfies $c\in C([0,T);L^\infty(\Omega))\cap L^\infty(Q_T)$. Furthermore, let $a\in C([0,T);W^{1,\infty}(\Omega))$ with $a,\nabla a$ uniformly bounded a.e. in $Q_T$. In other words, $a(\cdotp,\cdotp)$ satisfies the assumptions of the Picard--Lindel\"{o}f theorem: continuity with respect to $x,t$ and uniform Lipschitz continuity with respect to $x$.

\section{Exponential scaling}
\label{sec:exp_scaling}
The standard assumption on the coefficients $a,c$ found throughout the numerical literature is
\begin{equation}
c-\tfrac{1}{2}\diver a\ge\gamma_0>0 \quad\text{ on }Q_T 
\label{ellipt_assumpt_orig}
\end{equation}
for some constant $\gamma_0$. This assumption comes from the requirement that the weak formulations of the advection and reaction terms give an elliptic bilinear form on the corresponding function space. In \cite{Svadlenka}, this assumption is avoided by using \emph{exponential scaling} in time, i.e. the transformation $u(x,t)=e^{\alpha t}w(x,t)$ for $\alpha\in\RR$. After substitution into (\ref{contprob1}) and division of the whole equation by the common positive factor $e^{\alpha t}$ leads to a modified reaction coefficient $\tilde c=c+\alpha$ in the new equation for $w$. By choosing the constant $\alpha$ sufficiently large, (\ref{ellipt_assumpt_orig}) will be satisfied for the new equation. In \cite{Svadlenka}, error estimates that grow linearly in time are then derived for the DG scheme. However, for $\alpha>0$, if one transforms the resulting estimates back to the original problem using the exponential scaling transformation, the result is an estimate that depends exponentially on $T$.

Another possibility to avoid assumption (\ref{ellipt_assumpt_orig}) is a transformation similar to exponential scaling, however with respect to the spatial variable, cf. \cite{navert}, \cite{Johnson} and \cite{roos}: Let $\mu_0\in\RR^d$, then  we write
\begin{equation}
u(x,t)=e^{\mu_0\cdotp x}\tilde u(x,t).
\label{expon_scaling_const}
\end{equation}
Substituting into (\ref{contprob1}) and dividing the whole equation by the strictly positive function $e^{\mu_0\cdotp x}$, we get the new problem
\begin{equation}
\frac{\pd \tu}{\pd t}+a\cdotp\grad\, \tu+(a\cdotp\mu_0+c)\tu=0.
\label{contprob_scaled_const}
\end{equation}
The condition corresponding to (\ref{ellipt_assumpt_orig}) now reads: There exists $\mu_0\in\RR^d$ such that
\begin{equation}
a\cdotp\mu_0+c-\tfrac{1}{2}\diver a\ge\gamma_0>0 \quad\text{ on }Q_T.
\label{ellipt_assumpt_const}
\end{equation}

A possible generalization of the transformation (\ref{expon_scaling_const}) is using a function $\mu:\Omega\to\RR$ and setting
\begin{equation}
u(x,t)=e^{\mu(x)}\tilde u(x,t).
\label{expon_scaling_a}
\end{equation}
Again, substituting into (\ref{contprob1}) and dividing by $e^{\mu(x)}$, we get the new problem
\begin{equation}
\frac{\pd \tu}{\pd t}+a\cdot\grad \tu+(a\cdotp\grad\mu+c)\tu=0.
\label{contprob_scaled_a}
\end{equation}
The condition corresponding to (\ref{ellipt_assumpt_orig}) and (\ref{ellipt_assumpt_const}) now reads: There exists $\mu:\Omega\to\RR$ such that
\begin{equation}
a\cdotp\grad\mu+c-\tfrac{1}{2}\diver a\ge\gamma_0>0 \quad\text{ on }Q_T.
\label{ellipt_assumpt_a}
\end{equation}
This is essentially the approach used in \cite{Marini} for a stationary advection-diffusion-reaction problem. As shown in \cite{Devinatz}, the existence of a function $\mu:\Omega\to\RR$ such that $a\cdotp\grad\mu\ge \gamma_0$ is equivalent to the property that the advective field $a$ possesses neither closed curves nor stationary points. The uniformly positive term $a\cdotp\grad\mu$ can then be used to dominate the possibly negative term $c-\tfrac{1}{2}\diver a$ in order to satisfy condition (\ref{ellipt_assumpt_a}). If we used the choice (\ref{expon_scaling_a}) in our analysis, we would need to assume the nonexistence of closed curves or stationary points of the flow field for all $t$. This assumption is too restrictive.
 
In this paper, we shall generalize the transformation (\ref{expon_scaling_a}) using a function $\mu:Q_T\to\RR$ and set
\begin{equation}
u(x,t)=e^{\mu(x,t)}\tilde u(x,t).
\label{expon_scaling}
\end{equation}
Again, substituting into (\ref{contprob1}) and dividing by $e^{\mu(x,t)}$, we get the new problem
\begin{equation}
\frac{\pd \tu}{\pd t}+a\cdot\grad \tu+\Big(\frac{\pd \mu}{\pd t}+a\cdotp\grad\mu+c\Big)\tu=0.
\label{contprob_scaled}
\end{equation}
The condition corresponding to (\ref{ellipt_assumpt_orig}), (\ref{ellipt_assumpt_const}) and (\ref{ellipt_assumpt_a}) now reads: There exists $\mu:Q_T\to\RR$ such that
\begin{equation}
\frac{\pd \mu}{\pd t}+a\cdotp\grad\mu+c-\tfrac{1}{2}\diver a\ge\gamma_0>0 \quad\text{ on }Q_T.
\label{ellipt_assumpt}
\end{equation}
This is the condition that we will assume throughout the paper. Using the transformation (\ref{expon_scaling}) instead of (\ref{expon_scaling_a}) gives us sub-exponential growth or even uniform boundedness of the constant factor in the error estimate with respect to time.

If the original problem (\ref{contprob1}) does not satisfy (\ref{ellipt_assumpt_orig}), one is tempted to numerically solve the transformed equation (\ref{contprob_scaled}) instead, using DG or any other method. This is done e.g. in \cite{roos} for the case of linear $\mu$, i.e. (\ref{expon_scaling_const}). However then we are numerically solving a different problem and obtain different results, as the DG solutions of (\ref{contprob1}) and (\ref{contprob_scaled}) are not related by the simple relation (\ref{expon_scaling}), unlike the exact solutions. However as we will show in this paper, one can analyze the DG method for the original problem (\ref{contprob1}), while taking advantage of the weaker ellipticity condition (\ref{ellipt_assumpt}) for the transformed problem (\ref{contprob_scaled}).

Since the DG scheme is based on a suitable weak formulation, the first step is to reformulate the transformation (\ref{expon_scaling}) within the weak, rather than strong formulation of (\ref{contprob1}). The key step in deriving (\ref{contprob_scaled}) from (\ref{contprob1}) is dividing the whole equation by the common factor $e^{\mu(x,t)}$. However, if we substitute (\ref{expon_scaling}) into the weak formulation
\begin{equation}
\int_\Omega\frac{\pd u}{\pd t}v+a\cdotp\grad\, uv+cuv\dx=0,
\label{contprob_weak}
\end{equation}
it is not easy to divide the equation by $e^{\mu(x,t)}$, since this is inside the integral. The solution is to test (\ref{contprob_weak}) by the new test function  $\hat{v}(x,t)=e^{-\mu(x,t)}v(x,t)$. Due to the opposite signs in the exponents, the exponential factors cancel each other and we obtain the weak formulation of (\ref{contprob_scaled}):
\begin{equation}
\int_\Omega\frac{\pd \tu}{\pd t}v+a\cdotp\grad\, \tu v+\Big(\frac{\pd \mu}{\pd t} +a\cdotp\grad\mu+c\Big)\tu v\dx=0.
\label{contprob_weak_scaled}
\end{equation}
We note that the transformations $u\mapsto\tu$ and $v\mapsto\hat{v}$ are bijections. Furthermore, if we assume that the factor $e^{\mu(\cdot,t)}\in W^{1,\infty}(\Omega)$ as we shall do throughout the paper, the transformation (\ref{expon_scaling}) is a bijection from the function space $V$ into itself, where
\begin{equation}
V=\{u\in H^1(\Omega); u|_{\pdO^{-}}=0\}
\nonumber
\end{equation}
is the appropriate space for weak solutions and test functions for (\ref{contprob1}). This is the case $u_D=0$, the general case can be treated by the standard Dirichlet lifting procedure.

In the DG method, the above procedure using transform (\ref{expon_scaling}) cannot be directly applied, since if $v\in S_h$ -- the discrete space -- then  $\hat{v}=e^{-\mu}v$ will no longer lie in $S_h$ and therefore cannot be used as a test function in the formulation of the method. The solution is to test with a suitable projection of $\hat{v}$ onto $S_h$ and estimate the difference, as we shall do in Section \ref{sec:Anal_adv_react_terms}.

\section{Discontinuous Galerkin method}
\label{sec:DG}
Let $\Th$ be a triangulation of $\Omega$, i.e. a partition of $\overline{\Omega}$ into a finite number of closed simplices with mutually disjoint interiors. As usual with DG, $\Th$ need not be conforming, i.e. hanging nodes are allowed. For $K\in\Th$ we set $h_K=\mbox{diam}(K),\,\, h=\mbox{max}_{K\in\Th}h_K$.

For each $K\in\Th$ we define its \emph{inflow} and \emph{outflow} boundary by
\begin{equation}
\begin{split}
\partial K^-(t)&=\{x\in\partial K; a(x,t)\cdotp\bn(x)<0\},\\
\partial K^+(t)&=\{x\in\partial K; a(x,t)\cdotp\bn(x)\ge 0\},
\nonumber
\end{split}
\end{equation}
where $\bn(x)$ is the unit outer normal to $\partial K$. For simplicity, we shall omit the argument $t$ in the following and write simply $\partial K^\pm$.

Let $p\in\mathbb{N}\,$. The approximate solution will be sought in the space of discontinuous piecewise polynomial functions
\begin{equation}
S_h = \{v_h;\, v_h|_K \in P^p(K), \forall K\in\Th\},
\nonumber
\end{equation}
where $P^p(K)$ denotes the space of all polynomials on $K$ of degree at most $p$. Given an element $K\in\Th$, for $v_h\in S_h$ we define $v_h^-$ as the trace of $v_h$ on $\pd K$ from the side of the element adjacent to $K$. Furthermore on $\pd K\setminus\pd\Omega$ we define the \emph{jump} of $v_h$ as $[v_h]=v_h-v_h^-$, where $v_h$ is the trace from inside $K$.  

The DG formulation of (\ref{contprob1}) then reads: We seek $u_h\in C^1([0,T);S_h)$ such that $u_h(0)=u_h^0$, an $S_h$--approximation of the initial condition $u^0$, and for all $t\in(0,T)$
\begin{equation}
\Big(\frac{\partial u_h}{\partial t},v_h\Big)+b_h(u_h,v_h) +c_h(u_h,v_h)=l_h(v_h),\quad\forall v_h\in S_h, 
\label{discprob}
\end{equation}
where $b_h$ is the bilinear \emph{advection form}
\begin{equation}
b_h(u,v)= \sumK\intK (a\cdotp \grad u) v\dx -\sumK \intdKa (a\cdotp \bn)[u]v\dS -\sumK \intdKb (a\cdotp \bn)uv\dS,
\nonumber
\end{equation}
the \emph{reaction form} is defined by
\begin{equation}
c_h(u,v)= \intO cu v\dx
\nonumber
\end{equation}
and $l_h$ is the \emph{right-hand side form}
\begin{equation}
l_h(v)=  -\sumK \intdKb (a\cdotp \bn)u_Dv\dx.
\nonumber
\end{equation}
The definition of $b_h$ and $l_h$ corresponds to the concept of upwinding, cf. e.g. \cite{Svadlenka} for their derivation. 

\subsection{Auxiliary results}
In the following analysis we assume that the exact solution $u$ is sufficiently regular, namely
\begin{equation}
\label{ureg}
u,u_t:=\frac{\pd u}{\pd t}\in L^2(H^{p+1})
\end{equation}
We consider a system $\{\Th\}_{h\in(0,h_0)}$, $h_0>0$, of triangulations of ${\Omega}$ that are \emph{shape-regular} and satisfy the \emph{inverse assumption}, cf. \cite{CiarletNew}. Under these assumptions we have the following standard results:

\begin{lemma}[Inverse inequality]
\label{lem:inverseineq}
There exists a constant $C_I>0$ independent of $h,K$ such that for all $K\in\Th$, and all $v\in P^p(K)$,
\begin{equation}
|v|_{H^1(K)}\le C_I h_K^{-1} \|v\|_\LK.
\nonumber
\end{equation}
\end{lemma}


For $v\in L^{2}(\Omega)$ we denote by $\Pi_hv\in S_h$ the
$L^{2}(\Omega)$-projection of $v$ onto $S_h$:
\begin{equation}
\left(\Pi_h v-v,\,\varphi_h\right)=0,
    \quad \forall\,\varphi_h\in S_h.
\label{OG_def}
\end{equation}
Let $\eta_h(t)=u(t)-\Pi_hu(t)$ and $\xi_h(t)=\Pi_hu(t)-u_h(t)\in S_h$. Then we can write the error of the method as $e_h(t):=u(t)-u_h(t)=\eta_h(t_n)+\xi_h(t)$. For simplicity, we shall usually drop the subscript $h$. We have the following standard approximation result, cf. \cite{CiarletNew}.

\begin{lemma}
\label{lem4}
There exists a constant $C>0$ independent of $h,K$ such that for all $h\in (0,h_0)$ 
\begin{equation}
\begin{split}
\|\eta(t)\|&\le Ch^{p+1}|u(t)|_\Hppj,\\
|\eta(t)|_{H^1(K)}&\le Ch^{p}|u(t)|_\Hppj,\\
\|\eta(t)\|_{L^2(\pd K)}&\le Ch^{p+1/2}|u(t)|_\Hppj\\
\Big\|\frac{\partial\eta(t)}{\partial t}\Big\|&\leq Ch^{p+1}|u_t(t)|_\Hppj.
\nonumber
\end{split}
\end{equation}
\end{lemma}

\section{Analysis of the advection and reaction terms}
\label{sec:Anal_adv_react_terms}
In our analysis we will assume that there exists a constant $\gamma_0$ and a function $\mu:Q_T\to \RR$ such that (\ref{ellipt_assumpt}) holds. Furthermore, we assume that 
\begin{equation}
\begin{split}
0\leq\mu(x,t)&\leq\mu_{\max},\\
|\mu(x,t)-\mu(y,t)|&\leq L_\mu|x-y|,
\label{mu_assumpt}
\end{split}
\end{equation}
for all $x,y\in\Omega$ and $t\in (0,T)$. In other words, $\mu$ is nonnegative, uniformly bounded and Lipschitz continuous in space, where the Lipschitz constant is uniformly bounded for all $t$. Since $\Omega$ is a Lipschitz domain, hence quasi-convex, this means that $\mu(t)\in W^{1,\infty}(\Omega)$ for all $t$, with $W^{1,\infty}(\Omega)$ semi-norm uniformly bounded by $L_\mu$ for all $t$. We will show how to construct such a function in Section \ref{sec:mu_construct}.

Similarly as in Section \ref{sec:contprob}, we wish to write $\xi(x,t)=\emx\txi(x,t)$ and test the error equation with $\phi(x,t)=\emmx\txi(x,t)=\emdmx\xi(x,t)$ to obtain estimates for $\txi$. However, since $\phi(t)\notin S_h$ this is not possible. One possibility is to test by $\Pi_h\phi(t)\in S_h$ and estimate the resulting difference $\Pi_h\phi(t)-\phi(t)$. This is done in the stationary case in \cite{Marini} and the analysis is carried out under the assumption $\mu\in W^{p+1,\infty}(\Omega)$. Such high regularity can be achieved by mollification of $\mu$. However this would be somewhat technical in the evolutionary case, as space-time smoothing would be required in which case the dependance of all constants on $T$ must be carefully considered. Also $Q_T$ is potentially an unbounded domain (for $T=+\infty$) which leads to technical difficulties. Here we carry out the analysis under the weaker assumption (\ref{mu_assumpt}), i.e. $\mu(t)\in W^{1,\infty}(\Omega)$.

\begin{lemma}
\label{lem:phi_approx}
Let $\mu$ satisfy assumptions (\ref{mu_assumpt}). Let $\phi(x,t)=\emmx\txi(x,t)=\emdmx\xi(x,t)$, where $\xi(t)\in S_h$. Then there exists $C$ independent of $h,t,\xi,\txi$ such that
\begin{equation}
\begin{split}
\|\Pi_h\phi(t)-\phi(t)\|_\LK&\leq Ch_K \max_{x\in K}\emmx \|\txi(t)\|_\LK,\\
\|\Pi_h\phi(t)-\phi(t)\|_\LdK&\leq Ch^{1/2}_K \max_{x\in K}\emmx \|\txi(t)\|_\LK.
\label{lem:phi_approx_est}
\end{split}
\end{equation}
\end{lemma}

\proof
Let $x_K$ is the centroid of $K$. On element $K$ we introduce the constant $\mu_K(t)=\mu(x_K,t)$, then the function $e^{-2\mu_K(t)}\xi(\cdot,t)$ lies in $P^p(K)$, hence is fixed by the projection $\Pi_h$. Therefore,
\begin{equation}
\begin{split}
\Pi_h\phi(t)-\phi(t)&= \Pi_h\big(e^{-2\mu(t)}\xi(t) -e^{-2\mu_K(t)}\xi(t)\big) -\big(e^{-2\mu(t)}\xi(t) -e^{-2\mu_K(t)}\xi(t)\big)\\
&=\Pi_h w(t)-w(t),
\label{lem:phi_approx_est0}
\nonumber
\end{split}
\end{equation}
where $w(x,t)=(e^{-2\mu(x,t)} -e^{-2\mu(x_K,t)})\xi(x,t)$. Standard estimates of the interpolation error of $\Pi_h$ give
\begin{equation}
\|\Pi_h\phi(t)-\phi(t)\|_\LK^2 =\|\Pi_h w(t)-w(t)\|_\LK^2\leq Ch_K^{2}|w(t)|_{H^1(K)}^2.
\label{lem:phi_approx_est1}
\end{equation}
For the right-hand side seminorm we have
\begin{equation}
\begin{split}
|w(t)|_{H^1(K)}^2&=\int_K\big|\grad\big((e^{-2\mu(x,t)} -e^{-2\mu(x_K,t)})\xi(x,t)\big)\big|^2\dx\\
&\leq 2\int_K\big|\grad e^{-2\mu(x,t)}\xi(x,t)\big|^2\dx +2\int_K\big|(e^{-2\mu(x,t)} -e^{-2\mu(x_K,t)})\grad\xi(x,t)\big|^2\dx
\nonumber
\end{split}
\end{equation}
If $\mu(t)\in C^1(\overline{\Omega})$, by the mean value theorem
\begin{equation}
|\emdmxK-\emdmx|\leq h_K |\grad e^{-2\mu(\zeta,t)}| = h_K e^{-2\mu(\zeta,t)}2|\nabla\mu(\zeta,t)|
\nonumber
\end{equation}
for some point $\zeta$ on the line between $x$ and $x_K$. Therefore, by the inverse inequality
\begin{equation}
\begin{split}
|&w(t)|_{H^1(K)}^2\leq 2\int_K 4e^{-4\mu(x,t)} |\grad \mu(x,t)|^2|\xi(x,t)|^2\dx +2\int_K h_K^2e^{-4\mu(\zeta,t)}4|\nabla\mu(\zeta,t)|^2 |\grad\xi(x,t)|^2\dx\\
&\leq 8 |\mu(t)|^2_{W^{1,\infty}}\int_K e^{-4\mu(x,t)}e^{2\mu(x,t)} |\txi(x,t)|^2\dx +8C_I \max_{x\in K}e^{-4\mu(x,t)}|\mu(t)|^2_{W^{1,\infty}}  \|\xi(t)\|^2\\
&\leq 8 L_\mu^2 \max_{x\in K} e^{-2\mu(x,t)} e^{2h_KL_\mu} \|\txi(t)\|^2 +8 \max_{x\in K} e^{-2\mu(x,t)} e^{2h_KL_\mu} L_\mu^2 \|\txi(t)\|^2.
\label{lem:phi_approx_est2}
\end{split}
\end{equation}
Substituting into (\ref{lem:phi_approx_est1}) gives us the first inequality in (\ref{lem:phi_approx_est}) for $\mu(t)\in C^1(\overline{\Omega})$. The case $\mu(t)\in W^{1,\infty}(\Omega)$ follows by a density argument. The second inequality in (\ref{lem:phi_approx_est}) can be obtained similarly, only intermediately applying the trace inequality $\|\xi(t)\|_\LdK\leq Ch_K^{-1/2}\|\xi(t)\|_\LK$.
\qed

\begin{remark}
In Lemma \ref{lem:phi_approx}, the constant $C$ is proportional to $L_\mu e^{h_KL_\mu}$, as can be seen from (\ref{lem:phi_approx_est2}). Therefore the dependance of $C$ on $L_\mu$ is rather mild, effectively linear, due to the factor $h_K$ in the exponent.
\end{remark}

\medskip
Now we shall estimate individual terms in the DG formulation. Due to the consistency of the DG scheme, the exact solution $u$ also satisfies (\ref{discprob}). We subtract the formulations for $u$ and $u_h$ to obtain the error equation
\begin{equation}
\begin{split}
\Big(\frac{\partial \xi}{\partial t},v_h\Big) +\Big(\frac{\partial \eta}{\partial t},v_h\Big) + b_h(\xi,v_h) +b_h(\eta,v_h) + c_h(\xi,v_h) +c_h(\eta,v_h) =0
\label{error_eq}
\end{split}
\end{equation}
for all $v_h\in S_h$. 

As stated earlier we want to test (\ref{error_eq}) by $\phi(x,t)=\emmx\txi(x,t)$, however $\phi(t)\notin S_h$. We therefore set $v_h=\tPi_h \phi(t)$ and estimate the difference using Lemma \ref{lem:phi_approx}. We write (\ref{error_eq}) as
\begin{equation}
\begin{split}
\Big(\frac{\partial \xi}{\partial t},\Pi_h\phi\Big) &+b_h(\xi,\phi) +b_h(\xi,\Pi_h\phi-\phi)
+b_h(\eta,\phi) +b_h(\eta,\Pi_h\phi-\phi)\\ 
&+c_h(\xi,\phi) +c_h(\xi,\Pi_h\phi-\phi)
+c_h(\eta,\Pi_h\phi) +\Big(\frac{\partial \eta}{\partial t},\Pi_h\phi\Big)=0.
\label{error_eq_tested}
\end{split}
\end{equation}
We will estimate the individual terms of (\ref{error_eq_tested}) in a series of lemmas. For this purpose, we introduce the following norm on a subset $\omega$ of $\pd K$ or $\pd\Omega$:
\begin{equation}
\|f\|_{a,\omega}=\|\sqrt{|a\cdotp\bn|}f\|_{L^2(\omega)},
\nonumber
\end{equation}
where $\bn$ is the outer normal to $\pd K$ or $\pd\Omega$. We will usually omit the argument $t$ to simplify the notation.

\begin{lemma}
\label{lem:est2}
Let $\xi=e^\mu\txi, \phi=e^{-\mu}\txi$ as above and let $\mu$ satisfy assumption (\ref{ellipt_assumpt}). Then
\begin{equation}
\begin{split}
\Big(\frac{\partial \xi}{\partial t},\Pi_h\phi&\Big) +b_h(\xi,\phi)+c_h(\xi,\phi)\\
&\ge\frac{1}{2}\frac{\d}{\d t}\|\txi\|^2+\gamma_0\|\txi\|^2 +\frac{1}{2}\sumK\big(\big\|[\txi]\big\|_{a,\pd K^-\setminus\pd\Omega}^2 +\|\txi \|_{a,\pd K\cap\pd\Omega}^2\big).
\label{lem:est2:2}
\end{split}
\end{equation}
\end{lemma}

\proof
Since $\partial \xi/\partial t\in S_h$ for each $t$, by the  definition (\ref{OG_def}) of $\Pi_h$ we have
\begin{equation}
\Big(\frac{\partial \xi}{\partial t},\Pi_h\phi\Big)= \Big(\frac{\partial \xi}{\partial t},\phi\Big) =\Big(e^\mu\frac{\partial \txi}{\partial t} +e^\mu\frac{\partial \mu}{\partial t}\txi,e^{-\mu} \txi\Big) =\frac{1}{2}\frac{\d}{\d t}\|\txi\|^2 +\Big(\frac{\partial \mu}{\partial t}\txi,\txi\Big).
\label{lem:est2:1aaa}
\end{equation}
The reactive term satisfies
\begin{equation}
c_h(\xi,\phi)=\int_\Omega c\xi\phi\dx=\int_\Omega ce^\mu\txi e^{-\mu}\txi\dx =\int_\Omega c\txi^2\dx.
\label{lem:est2:1aa}
\end{equation}
From the definition of $b_h$, we get
\begin{equation}
\begin{split}
b_h(\xi,\phi)&=\!\! \sumK\intK a\cdotp(\grad\mu\,\txi+\grad\txi)e^\mu e^{-\mu}\txi\dx\\ 
&\quad-\!\!\sumK \intdKa (a\cdotp \bn)[e^\mu\txi]e^{-\mu}\txi\dS -\!\!\sumK \intdKb (a\cdotp \bn)e^\mu\txi e^{-\mu}\txi\dS\\
&=\!\! \sumK\intK a\cdotp(\grad\mu\,\txi+\grad\txi)\txi\dx -\!\!\sumK \intdKa (a\cdotp \bn)[\txi]\txi\dS -\!\!\sumK \intdKb (a\cdotp \bn)\txi^2\dS.
\label{lem:est2:1a}
\end{split}
\end{equation}
By Green's theorem,
\begin{equation}
\intK a\cdotp\!\grad\txi\,\txi\dx=-\frac{1}{2}\int_K \diver a\,\txi^2\dx +\frac{1}{2}\int_{\pd K}(a\cdotp\bn)\,\txi^2\dS.
\label{lem:est2:1b}
\end{equation}
Splitting the last integral over the separate parts of $\pd K$: $\pd K^-\!\setminus\!\pd\Omega, \pd K^-\!\cap\Omega, \pd K^+\!\setminus\!\pd\Omega$, and $\pd K^+\!\cap\Omega$, and by substituting (\ref{lem:est2:1b}) into $(\ref{lem:est2:1a})$, we get 
\begin{equation}
\begin{split}
b_h(\xi,\phi)&=\!\! \sumK\intK (a\cdotp\grad\mu-\tfrac{1}{2}\diver a)\txi^2\dx\\ &+\!\!\sumK\bigg(-\frac{1}{2}\intdKa (a\cdotp \bn)(\txi^2-2\txi\txi^-)\dS -\frac{1}{2}\intdKb (a\cdotp \bn)\txi^2\dS\\ 
&\qquad+\frac{1}{2}\int_{\pd K^+\!\setminus\pd\Omega}(a\cdotp \bn)\txi^2\dS +\frac{1}{2}\int_{\pd K^+\!\cap\pd\Omega}(a\cdotp \bn)\txi^2\dS\bigg).
\label{lem:est2:1c}
\end{split}
\end{equation}
We note that
\begin{equation}
\sumK\int_{\pd K^+\!\setminus\pd\Omega}(a\cdotp \bn)\txi^2\dS =-\sumK\int_{\pd K^-\!\setminus\pd\Omega}(a\cdotp \bn)(\txi^-)^2\dS.
\label{lem:est2:1d}
\end{equation}
Using the facts that $\txi^2-2\txi\txi^-+(\txi^-)^2=[\txi]^2$ and   $-a\cdotp\bn=|a\cdotp\bn|$ on $\pd K^-$ and $a\cdotp\bn=|a\cdotp\bn|$ on $\pd K^+$, we get (\ref{lem:est2:2}) by substituting (\ref{lem:est2:1d}) into (\ref{lem:est2:1c}) and applying assumption (\ref{ellipt_assumpt}) in the resulting interior terms along with (\ref{lem:est2:1aaa}), (\ref{lem:est2:1aa}).

\qed

\begin{lemma}
\label{lem:est3}
We have
\begin{equation}
|b_h(\xi,\Pi_h\phi-\phi)|\leq Ch\|\txi\|^2+\frac{1}{8}\sumK\big(\big\|[\txi]\big\|_{a,\pd K^-\!\setminus\pd\Omega}^2 +\|\txi \|_{a,\pd K\cap\pd\Omega}^2\big).
\nonumber
\end{equation}
\end{lemma}
\proof
We estimate the terms of $b_h(\xi,\Pi_h\phi-\phi)$ over the interiors and boundaries of elements separately. Let $\Pi_h^1$ be the $L^2(\Omega)$-projection onto the space of discontinuous piecewise linear polynomials on $\Th$. Since on each $K\in\Th$ it holds that $\grad\txi|_K\in P^{p-1}(K)$, then $\Pi_h^1 a\cdotp\grad\txi\in S_h$. Hence due to (\ref{OG_def}), we have
\begin{equation}
\sum_K\int_K\Pi_h^1 a\cdotp\grad\txi(\Pi_h\phi-\phi)\dx=0.
\nonumber
\end{equation}
Due to standard approximation results, we have $\|a-\Pi_h^1 a\|_{L^\infty(K)}\leq Ch_K|a|_{W^{1,\infty}(K)}$, thus we can estimate the interior terms of $b_h(\xi,\Pi_h\phi-\phi)$ as
\begin{equation}
\begin{split}
\sumK&\intK a\cdotp\grad\xi(\Pi_h\phi-\phi)\dx =\sumK\intK (a-\Pi_h^1 a)\cdotp\grad\xi(\Pi_h\phi-\phi)\dx\\
&\leq\sumK Ch_K|a|_{W^{1,\infty}}C_Ih_K^{-1}\|\xi\|_{L^2(K)} \|\Pi_h\phi-\phi\|_{L^2(K)}\\
&\leq C\sumK \max_{x\in K}\emx \|\txi\|_{L^2(K)} Ch \max_{x\in K}\emmx\|\txi\|_\LK\leq Ce^{L_\mu h}h\|\txi\|^2 \leq Ch\|\txi\|^2,
\label{lem:est3:1a}
\end{split}
\end{equation}
due to the inverse inequality and Lemma \ref{lem:phi_approx}. For the boundary terms of $b_h(\xi,\Pi_h\phi-\phi)$, we get
\begin{equation}
\begin{split}
-&\sumK \intdKa (a\cdotp \bn)[\xi](\Pi_h\phi-\phi)\dS -\sumK \intdKb (a\cdotp \bn)\xi(\Pi_h\phi-\phi)\dS\\ &\leq\sumK \intdKa |a\cdotp \bn|\max_{x\in K}(\emx)\big|[\txi]\big||\Pi_h\phi-\phi| \dS +\!\!\sumK \intdKb |a\cdotp \bn|\max_{x\in K}(\emx)|\txi||\Pi_h\phi-\phi|\dS\\
&\leq \frac{1}{8}\sumK \intdKa |a\cdotp \bn|[\txi]^2\dS +\frac{1}{8}\sumK \intdKb |a\cdotp \bn|\txi^2\dS\\
&\quad+C\sumK\max_{x\in K}\edmx\intdK|a\cdotp\bn||\Pi_h\phi-\phi|^2\dS\\
&\leq \frac{1}{8}\sumK\big(\big\|[\txi]\big\|_{a,\pd K^-\!\setminus\pd\Omega}^2 +\|\txi \|_{a,\pd K\cap\pd\Omega}^2\big) + C\sum_K\max_{x\in K}\edmx h_K \max_{x\in K}\emdmx\|\txi\|_\LK^2
\label{lem:est3:1b}
\end{split}
\end{equation}
by Young's inequality and Lemma \ref{lem:phi_approx}. Again we estimate $\max_{x\in K}\edmx\max_{x\in K}\emdmx\leq e^{2L_\mu h_K}$ in (\ref{lem:est3:1b}), which completes the proof after combining with (\ref{lem:est3:1a}).

\qed

\begin{lemma}
\label{lem:est4}
We have
\begin{equation}
|b_h(\eta,\phi)|\leq Ch^{p+1}|u(t)|_\Hppj\|\txi\| + Ch^{2p+1}|u(t)|_\Hppj^2 +\frac{1}{8}\sumK\big(\big\|[\txi]\big\|_{a,\pd K^-\!\setminus\pd\Omega}^2 +\|\txi \|_{a,\pd K\cap\pd\Omega}^2\big).
\nonumber
\end{equation}
\end{lemma}

\proof
We have by Green's theorem
\begin{equation}
\begin{split}
b_h(\eta,\phi)= \sumK&\bigg(\intdK a\cdotp\bn\eta\phi\dS -\intK (\diver a)\eta\phi\dx-\intK a\cdotp\grad\phi\eta\dx\\
&-\intdKa (a\cdotp \bn)[\eta]\phi\dS -\intdKb (a\cdotp \bn)\eta\phi\dS\bigg).
\label{lem:est4:1a}
\end{split}
\end{equation}
The first integral over $K$ can be estimated as
\begin{equation}
-\sumK\intK (\diver a)\eta\phi\dx\leq Ch^{p+1}|u(t)|_\Hppj\max_{x\in \Omega}\emmx\|\txi\| \leq Ch^{p+1}|u(t)|_\Hppj\|\txi\|,
\nonumber
\end{equation}
because $\mu\ge 0$. Since $\phi=e^{-2\mu}\xi$, we get for the second integral over $K$ in (\ref{lem:est4:1a}):
\begin{equation}
\begin{split}
-&\sumK\intK a\cdotp\grad\phi\,\eta\dx =\sumK\intK 2a\cdotp\grad\mu\,e^{-2\mu}\xi\eta\dx -\sumK\intK e^{-2\mu}a\cdotp\grad\xi\,\eta\dx.
\label{lem:est4:1b}
\end{split}
\end{equation}
The first right-hand side term in (\ref{lem:est4:1b}) can be estimated by
\begin{equation}
\begin{split}
\sumK&\intK 2a\cdotp\grad\mu\,e^{-2\mu}\xi\eta\dx =\sumK\intK 2a\cdotp\grad\mu\,e^{-\mu}\txi\eta\dx\\
&\leq\sumK CL_\mu\max_{x\in K}e^{-\mu(x,t)}h^{p+1}|u(t)|_{\Hppj(K)} \|\txi\|_{L^2(K)}
\leq CL_\mu h^{p+1}|u(t)|_{\Hppj}\|\txi\|.
\nonumber
\end{split}
\end{equation}
The second right-hand side term in (\ref{lem:est4:1b}) can be estimated similarly as in (\ref{lem:est3:1a}), due to the definition of $\eta$:
\begin{equation}
\begin{split}
-&\sumK\intK e^{-2\mu}a\cdotp\grad\xi\,\eta\dx =\sumK\intK \big(\Pi_h^1(e^{-2\mu}a)-e^{-2\mu}a\big) \cdotp\grad\xi\eta\dx\\
&\leq\sumK Ch_K|e^{-2\mu}a|_{W^{1,\infty}}C_Ih_K^{-1}\|\xi\|_{L^2(K)} h^{p+1}|u(t)|_{H^{p+1}(K)}\leq C\max_{x\in \Omega}e^{-3\mu(x,t)}h^{p+1}|u(t)|_{\Hppj}\|\txi\|\\
&\leq Ch^{p+1}|u(t)|_{\Hppj}\|\txi\|,
\nonumber
\end{split}
\end{equation}
since $|e^{-2\mu}a|_{W^{1,\infty}}=|-2\grad\mu\, e^{-2\mu}a +e^{-2\mu}\diver a|_{L^{\infty}}\leq C$.

As for the boundary terms in (\ref{lem:est4:1a}), we can split the integral over  $\pd K$ into integrals over the separate parts $\pd K^-\!\setminus\!\pd\Omega, \pd K^-\!\cap\pd\Omega, \pd K^+\!\setminus\!\pd\Omega$ and $\pd K^+\!\cap\pd\Omega$, similarly as in the proof of Lemma \ref{lem:est2}. Thus several terms are cancelled out:
\begin{equation}
\begin{split}
\sumK&\bigg(\intdK (a\cdotp\bn)\eta\phi\dS -\intdKa (a\cdotp \bn)[\eta]\phi\dS -\intdKb (a\cdotp \bn)\eta\phi\dS\bigg)\\
&=\sumK\bigg(\int_{\pd K^+\setminus\pd\Omega} (a\cdotp\bn)\eta\phi\dS +\int_{\pd K^+\cap\pd\Omega} (a\cdotp\bn)\eta\phi\dS +\intdKa (a\cdotp \bn)\eta^-\phi\dS\bigg)\\
&=\sumK\bigg(\int_{\pd K^-\setminus\pd\Omega} (a\cdotp\bn)\eta^-[\phi]\dS +\int_{\pd K^+\cap\pd\Omega} (a\cdotp\bn)\eta\phi\dS\bigg)
\label{lem:est4:1d}
\end{split}
\end{equation}
using a similar identity to (\ref{lem:est2:1d}). Finally, we can use Young's inequality to estimate (\ref{lem:est4:1d}) further as
\begin{equation}
\begin{split}
\ldots&\leq\sumK\bigg(\int_{\pd K^-\setminus\pd\Omega} |a\cdotp\bn||\eta^-|e^{-\mu}\big|[\txi]\big|\dS +\int_{\pd K^+\cap\pd\Omega} |a\cdotp\bn||\eta|e^{-\mu}|\txi|\dS\bigg)\\
&\leq \frac{1}{8}\sumK\bigg(\int_{\pd K^-\setminus\pd\Omega} |a\cdotp\bn|[\txi]^2\dS +\int_{\pd K^+\cap\pd\Omega} |a\cdotp\bn|\txi^2\dS\bigg) +C\sum_K\intdK|a\cdotp\bn|\eta^2e^{-2\mu}\dS\\
&\leq \frac{1}{8}\sumK\big(\big\|[\txi]\big\|_{a,\pd K^-\!\setminus\pd\Omega}^2 +\|\txi \|_{a,\pd K\cap\pd\Omega}^2\big) +Ch^{2p+1}|u(t)|_\Hppj^2.
\nonumber
\end{split}
\end{equation}
The proof is completed by gathering all the above estimates of the individual terms of $b_h(\eta,\phi)$.

\qed

\begin{lemma}
\label{lem:est5}
We have
\begin{equation}
|b_h(\eta,\Pi_h\phi-\phi)|\leq Ch^{p+1}|u(t)|_\Hppj\|\txi\|.
\nonumber
\end{equation}
\end{lemma}

\proof
We use Lemmas \ref{lem4} and \ref{lem:phi_approx} to estimate
\begin{equation}
\begin{split}
b_h(\eta&,\Pi_h\phi-\phi) =\sumK\bigg(\intK (a\cdotp \grad \eta) (\Pi_h\phi-\phi)\dx\\
&\quad-\intdKa (a\cdotp \bn)[\eta](\Pi_h\phi-\phi)\dS -\intdKb (a\cdotp \bn)\eta (\Pi_h\phi-\phi)\dS\bigg)\\
&\leq Ch^p|u(t)|_\Hppj Ch \max_{x\in\Omega}\emmx\|\txi\| +Ch^{p+1/2}|u(t)|_\Hppj Ch^{1/2} \max_{x\in\Omega}\emmx\|\txi\|.
\nonumber
\end{split}
\end{equation}

\qed

\begin{lemma}
\label{lem:est6}
We have
\begin{equation}
\begin{split}
|c_h(\xi,\Pi_h\phi-\phi)|&\leq Ch\|\txi\|^2,\\ 
|c_h(\eta,\Pi_h\phi)|&\leq Ch^{p+1}|u(t)|_\Hppj\|\txi\|.
\nonumber
\end{split}
\end{equation}
\end{lemma}

\proof
Lemma \ref{lem:phi_approx} gives us
\begin{equation}
\sumK\int_K c\xi(\Pi_h\phi-\phi)\dx\leq C\sumK\max_{x\in K}e^{\mu(x,t)}\|\txi\|_{L^2(K)} h_K \max_{x\in K}\emmx\|\txi\|_\LK \leq Ce^{hL_\mu} h\|\txi\|^2.
\nonumber
\end{equation}
As for the second estimate, we write $c_h(\eta,\Pi_h\phi) =c_h(\eta,\phi)+c_h(\eta,\Pi_h\phi-\phi)$ and estimate by Lemmas \ref{lem4} and \ref{lem:phi_approx}:
\begin{equation}
\begin{split}
c_h(\eta,\phi) &=\int_\Omega c\eta\phi\dx=\int_\Omega \eta e^{-\mu}\txi\dx \leq  Ch^{p+1}|u(t)|_\Hppj\|\txi\|,\\
c_h(\eta,\Pi_h\phi-\phi) &=\intO c\eta (\Pi_h\phi-\phi)\dx
\leq Ch^{p+1}|u(t)|_\Hppj Ch \|\txi\|.
\nonumber
\end{split}
\end{equation}
Combining these two estimates gives the desired result.
\qed

\begin{lemma}
\label{lem:est7}
We have
\begin{equation}
\Big|\Big(\frac{\partial \eta}{\partial t},\Pi_h\phi\Big)\Big| \leq Ch^{p+1}|u_t(t)|_\Hppj\|\txi\|.
\nonumber
\end{equation}
\end{lemma}

\proof
We use Lemmas \ref{lem4}, \ref{lem:phi_approx} and the definition of $\phi$:
\begin{equation}
\begin{split}
\Big(\frac{\partial \eta}{\partial t}&,\Pi_h\phi\Big)= \Big(\frac{\partial \eta}{\partial t},\phi\Big) +\Big(\frac{\partial \eta}{\partial t},\Pi_h\phi-\phi\Big)
\leq Ch^{p+1}|u_t(t)|_\Hppj\|\txi\| +Ch^{p+1}|u_t(t)|_\Hppj Ch \|\txi\|.
\nonumber
\end{split}
\end{equation}
\qed

\section{Error analysis}
\label{sec:MethLines}
Finally we come to the error analysis. The starting point is the error identity (\ref{error_eq_tested}) to which we apply the derived estimates of its individual terms.

\begin{theorem}
\label{thm:main}
Let there exist a function $\mu:Q_T\to[0,\mu_{\max}]$ for some constant $\mu_{\max}$, such that $\mu(t)\in W^{1,\infty}(\Omega)$ and let there exist a constant $\gamma_0>0$ such that the coefficients of (\ref{contprob1}) satisfy  $\mu_t +a\cdotp\grad\mu+c-\tfrac{1}{2}\diver a\ge\gamma_0>0$ on $Q_T$. Let the initial condition $u_h^0$ satisfy $\|\Pi_h u^0-u_h^0\|\leq Ch^{p+1/2}|u^0|_\Hppj$. Then there exists a constant $C$ depending on $a,c,\mu$ but independent of $h$ and $T$ such that the error of the DG scheme (\ref{discprob}) satisfies
\begin{equation}
\begin{split}
\max_{t\in[0,T]}\|e_h(t)\| &+\sqrt{\gamma_0}\|e_h\|_{L^2(Q_T)} +\frac{1}{2}\Big(\int_0^T\sumK\big(\big\|[e_h(\vt)]\big\|_{a,\pd K^-\setminus\pd\Omega}^2 +\|e_h(\vt)\|_{a,\pd K\cap\pd\Omega}^2\big)\d\vt\Big)^{1/2}\\
&\leq Ch^{p+1/2}\big(|u^0|_\Hppj+|u|_{L^2(\Hppj)} +|u_t|_{L^2(\Hppj)}\big).
\label{thm:main:est}
\end{split}
\end{equation}
\end{theorem}

\proof
Applying Lemmas \ref{lem:est2} -- \ref{lem:est7}  to (\ref{error_eq_tested}), multiplying by 2 for convenience and collecting similar terms, we get for all $t$
\begin{equation}
\begin{split}
\frac{\d}{\d t}&\|\txi(t)\|^2 +2\gamma_0\|\txi(t)\|^2 +\frac{1}{2}\sumK\big(\big\|[\txi(t)]\big\|_{a,\pd K^-\setminus\pd\Omega}^2 +\|\txi(t)\|_{a,\pd K\cap\pd\Omega}^2\big)\\
&\leq Ch\|\txi(t)\|^2 +Ch^{2p+1}|u(t)|_\Hppj^2 +Ch^{p+1}\big(|u_t(t)|_\Hppj +|u(t)|_\Hppj\big)\|\txi(t)\|\\
&\leq Ch\|\txi(t)\|^2 +Ch^{2p+1}\big(|u(t)|_\Hppj^2+|u_t(t)|_\Hppj^2\big),
\nonumber
\end{split}
\end{equation}
where we have applied Young's inequality to the last right-hand side product on the second line. If $h\leq c_h$  for $c_h$ sufficiently small, we have $Ch\leq \gamma_0$, hence the first right-hand side term can be hidden under the left-hand side term $2\gamma_0\|\txi(t)\|^2$:
\begin{equation}
\begin{split}
\frac{\d}{\d t}\|\txi(t)\|^2 &+\gamma_0\|\txi(t)\|^2 +\frac{1}{2}\sumK\big(\big\|[\txi(t)]\big\|_{a,\pd K^-\setminus\pd\Omega}^2 +\|\txi(t)\|_{a,\pd K\cap\pd\Omega}^2\big)\\
&\leq Ch^{2p+1}\big(|u(t)|_\Hppj^2+|u_t(t)|_\Hppj^2\big).
\nonumber
\end{split}
\end{equation}
Substituting $t=\vt$ and integrating over $(0,t)$ gives us
\begin{equation}
\begin{split}
\|\txi(t)\|^2 &+\gamma_0\int_0^t\|\txi(\vt)\|^2\d \vt +\frac{1}{2}\int_0^t\sumK\big(\big\|[\txi(\vt)]\big\|_{a,\pd K^-\setminus\pd\Omega}^2 +\|\txi(\vt)\|_{a,\pd K\cap\pd\Omega}^2\big)\d\vt\\
&\leq \|\txi(0)\|^2 +Ch^{2p+1}\big(|u|_{L^2(0,t;\Hppj)}^2 +|u_t|_{L^2(0,t;\Hppj)}^2\big)\\ 
&\leq Ch^{2p+1}\big(|u^0|_\Hppj^2 +|u|_{L^2(0,t;\Hppj)}^2 +|u_t|_{L^2(0,t;\Hppj)}^2\big),
\label{thm:main_2}
\end{split}
\end{equation}
since the assumptions give
\begin{equation} 
\|\txi(0)\|^2\leq \|\xi(0)\|^2 =\|\Pi_h u^0-u_h^0\|^2 \leq Ch^{2p+1}|u^0|_\Hppj^2.
\nonumber
\end{equation}
Now we reformulate estimate (\ref{thm:main_2}) as an estimate of $\xi$ instead of $\txi$. Because $\txi=e^{-\mu}\xi$, we can estimate for example
\begin{equation}
\|\txi(t)\|^2\ge \min_{Q_T}e^{-2\mu(x,t)}\|\xi(t)\|^2 =e^{-2\max_{Q_T}\mu(x,t)}\|\xi(t)\|^2 
\nonumber
\end{equation}
and similarly for the remaining left-hand side norms in (\ref{thm:main_2}). Substituting into (\ref{thm:main_2}) and multiplying by $e^{2\max_{Q_T}\mu(x,t)}= e^{2\mu_{\max}}$ gives us
\begin{equation}
\begin{split}
\|\xi(t)\|^2 &+\gamma_0\|\xi(\vt)\|_{L^2(0,t;L^2)}^2 +\frac{1}{2}\int_0^t\sumK\big(\big\|[\xi(\vt)]\big\|_{a,\pd K^-\setminus\pd\Omega}^2 +\|\xi(\vt)\|_{a,\pd K\cap\pd\Omega}^2\big)\d\vt\\
&\leq Ch^{2p+1} e^{2\mu_{\max}}\big(|u^0|_\Hppj^2 +|u|_{L^2(0,t;\Hppj)}^2 +|u_t|_{L^2(0,t;\Hppj)}^2\big).
\label{thm:main_3}
\end{split}
\end{equation}
Lemma \ref{lem4} gives a similar inequality for $\eta$ and the triangle inequality gives us (\ref{thm:main:est}).
\qed

\section{Construction of the function $\mu$}
\label{sec:mu_construct}
In this section we show a construction of the function $\mu$ satisfying (\ref{ellipt_assumpt}). Many different constructions of $\mu$ are possible depending on the assumptions on the vector field $a$. For example we can always take $\mu(x,t)=\alpha t$ for some $\alpha\ge 0$. This corresponds to standard exponential scaling and, as we have seen, this choice leads to exponential growth in time of the error estimate. Another possibility was the approach of \cite{Marini} mentioned in Section \ref{sec:contprob}, where a suitable function $\mu(x)$ exists if $a$, which is stationary, possesses no closed curves or stationary points. Here we show another possibility with an interesting interpretation.

If $c-\tfrac{1}{2}\diver a$ is negative or changes sign frequently, we can use the expression $\mu_t+a\cdot\grad\mu$ to dominate this term everywhere. If we choose $\mu_1$ such that
\begin{equation}
\frac{\pd\mu_1}{\pd t}+a\cdot\grad\mu_1=1\quad \text{on }Q_T,
\label{mu_construct:1}
\end{equation}
then by multiplying $\mu_1$ by a sufficiently large constant, we can satisfy the ellipticity condition (\ref{ellipt_assumpt}) for a chosen $\gamma_0>0$.

Equation (\ref{mu_construct:1}) can be explicitly solved using characteristics. We define \emph{pathlines} of the flow, i.e. the family of curves $S(t;x_0,t_0)$ by
\begin{equation}
S(t_0;x_0,t_0)=x_0\in\overline{\Omega},\quad \frac{\d S(t;x_0,t_0)}{\d t} =a(S(t;x_0,t_0),t).
\label{mu_construct:2}
\end{equation}
This means that $S(\cdot;t_0,x_0)$ is the trajectory of a massless particle in the nonstationary flow field $a$ passing through point $x_0$ at time $t_0$. It is convenient to choose the parameter $t_0$  minimal for each pathline -- then the pair $(x_0,t_0)$ is the ``origin" of the pathline. In other words, for each $x_0\in\overline{\Omega}$, there is a pathline $S(\cdot;x_0,0)$ corresponding to trajectories of particles present in $\Omega$ at the initial time $t_0=0$. Then there are trajectories of particles entering through the inlet part of $\partial\Omega$: for each $x_0\in\pd\Omega^-$ and all $t_0$ there exists a pathline $S(\cdot;x_0,t_0)$ originating at $(x_0,t_0)\in\pd\Omega^-\times[0,T)$.

Equation (\ref{mu_construct:1}) can then be rewritten along the pathlines as
\begin{equation}
\frac{\d\, \mu_1(S(t;x_0,t_0),t)}{\d t} = \Big(\frac{\pd\mu_1}{\pd t}+a\cdot\grad\mu_1\Big) (S(t;x_0,t_0),t)=1,
\nonumber
\end{equation}
therefore
\begin{equation}
\mu_1(S(t;x_0,t_0),t) =t-t_0.
\label{mudef}
\end{equation}
Here any constant can be chosen instead of $t_0$, however this choice is the most convenient. Since then at the origin of a pathline we have
\begin{equation}
\mu_1(S(t_0;x_0,t_0),t_0) =0
\nonumber
\end{equation}
and the value of $\mu_1$ along this pathline is simply the time elapsed since $t_0$. In other words, $\mu_1(x,t)$ is the time a particle carried by the flow passing through $x\in\Omega$ at time $t$ has spent in $\Omega$ up to the time $t$. In this paper, we assume this quantity to be uniformly bounded in order to satisfy assumption (\ref{mu_assumpt}).

Now that we have constructed a function satisfying (\ref{mu_construct:1}), we can choose $\gamma_0$ and define e.g.
\begin{equation}
\mu(x,t)=\mu_1(x,t)\big(\big|\inf_{Q_T}(c-\tfrac{1}{2}\diver a)^-\big|+\gamma_0\big),
\nonumber
\end{equation}
where $f^-=\min\{0,f\}$ is the negative part of $f$. Then 
\begin{equation}
\frac{\pd\mu}{\pd t} +a\cdot\grad\mu+c-\tfrac{1}{2}\diver a = \big|\inf_{Q_T}(c-\tfrac{1}{2}\diver a)^-\big|+\gamma_0 +c-\tfrac{1}{2}\diver a \ge\gamma_0.
\nonumber
\end{equation}
This choice of $\mu$ leads to estimates of the following form, cf. (\ref{thm:main_3}),
\begin{equation}
\|e_h\|_{L^\infty(L^2)} +\sqrt{\gamma_0}\|e_h\|_{L^2(L^2)}\leq Ce^{\widehat{T}(|\inf_{Q_T}(c-\tfrac{1}{2}\diver a)^-|+\gamma_0)}h^{p+1/2},
\label{sec:mu_construct:est1}
\end{equation}
where $\widehat{T}=\sup_{Q_T}\mu_1$, i.e. $\widehat{T}$ is the maximal time any particle carried by the flow field $a$ spends in $\Omega$. If we compare this to estimates obtained by a straightforward analysis using Gronwall's lemma without any ellipticity assumption, we would expect
\begin{equation}
\|e_h\|_{L^\infty(L^2)} \leq Ce^{T(\sup_{Q_T}|c-\tfrac{1}{2}\diver a|)}h^{p+1/2}.
\label{sec:mu_construct:est2}
\end{equation}
Comparing (\ref{sec:mu_construct:est2}) to (\ref{sec:mu_construct:est1}), we see that the estimates are essentially of similar form, only the exponential dependence on global physical time $T$ has been replaced by dependence on time $\widehat{T}$ along pathlines, which is bounded. Effectively, our analysis replaces the application of Gronwall's lemma in the Eulerian framework with its application in the Lagrangian framework -- along individual pathlines which have bounded length.

Throughout the paper, we have assumed $\widehat{T}$ to be bounded. In general, we could consider dependencies of the form  
\begin{equation}
\widehat{T}(t)=\sup_{(x,\vartheta)\in\Omega\times{(0,t)}} \mu_1(x,t),
\nonumber
\end{equation}
i.e. $\widehat{T}(t)$ is the maximal time any particle carried by the flow $a$ spends in $\Omega$ up to time $t$. From (\ref{sec:mu_construct:est1}), we can expect exponential dependence of the $L^\infty(0,t;L^2)$ and $L^2(0,t;L^2)$ norms on $\widehat{T}(t)$. The result of the standard analysis (\ref{sec:mu_construct:est2}) corresponds to the ``worst case" $\widehat{T}(t)=t$, i.e. that there is a particle that stays inside $\Omega$ for all $t\in[0,T)$. However considering more general dependencies on time is possible, e.g. $\widehat{T}=\sqrt{t}$, leading to growth of the error which is exponential in $\sqrt{t}$.

\subsection{Regularity of the function $\mu$} 
\begin{figure}[t!]
\centering\includegraphics[scale=1.,trim={0 0 3.cm 0}, clip]{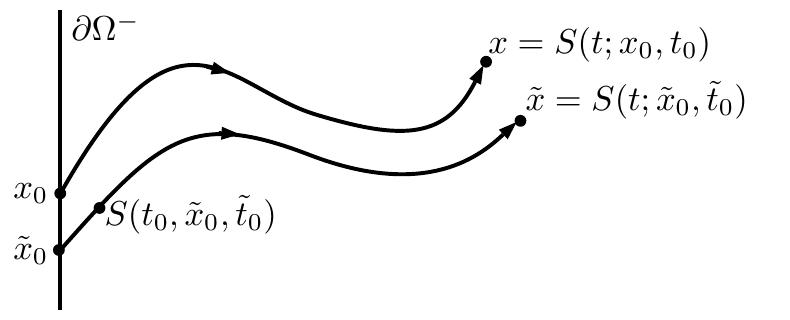}
\includegraphics[scale=1.,trim={0 0 1.cm 0}, clip]{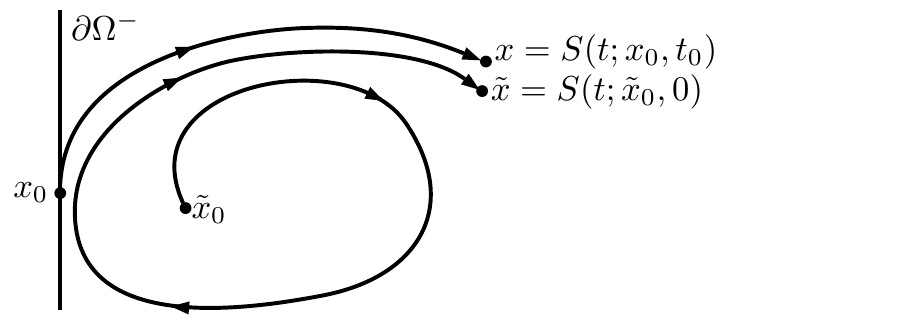}
\caption{Left: proof of Lemma \ref{lem:mu_reg}. Right: vortex touching $\pd\Omega^-$ with $t_0\gg\tt_0=0$, i.e. $t_0-\tt_0$ is large, hence $\mu_1$ is not Lipschitz continuous at $x$.}
\label{fig:1}
\end{figure}

Now we show when $\mu_1$ defined by (\ref{mudef}) satisfies conditions (\ref{mu_assumpt}), especially Lipschitz continuity in space which was necessary in the analysis. One obvious example when $\mu_1$ defined by (\ref{mudef}) will not be Lipschitz continuous is when a vortex will ``touch" $\pd\Omega^-$ as in Figure \ref{fig:1}. Then if $(x,t)$ lies on a pathline originating at $(x_0,t_0)$ such that $a(x_0,t_0)$ is tangent to $\pd\Omega$, then we can find $\tx$ arbitrarily close to $x$ such that the corresponding pathline is much longer, perhaps winding several times around the vortex. Then $\mu_1(x,t)-\mu_1(\tx,t)=\tt_0-t_0$ can be very large, while $\|x-\tx\|$ is arbitrarily small, hence $\mu_1$ is not Lipschitz continuous. In fact $\mu_1$ is discontinuous at $(x,t)$. In the following lemma, we show that inlet points where $a$ is tangent to $\pd\Omega$ are the only troublemakers.

We note that since $\Omega$ is a Lipschitz domain (hence quasiconvex), when proving Lipschitz continuity of $\mu_1$ in space it is sufficient to prove local Lipschitz continuity of $\mu_1$ in some neighborhood of each $x\in\overline{\Omega}$ with a Lipschitz constant independent of $x$.

\begin{lemma}
\label{lem:mu_reg}
Let $a\in L^{\infty}(Q_T)$ be continuous with respect to time and Lipschitz continuous with respect to space. Let there exist a constant $a_{\min}>0$ such that 
\begin{equation}
-a(x,t)\cdot\bn\ge a_{\min}
\label{lem:mu_reg_ass}
\end{equation}
for all $x\in\pd\Omega^-,t\in[0,T)$. Let $\mu_1$ be defined by (\ref{mudef}) on $\overline{\Omega}\times[0,T)$. Let the time any particle carried by the flow field $a(\cdot,\cdot)$ spends in $\Omega$ be uniformly bounded by $\widehat{T}$. Then $\mu_1$ satisfies assumption (\ref{mu_assumpt}).
\end{lemma}
\proof
By definition, we have $\mu(x,t)\ge 0$ for all $x,t$. By the above considerations, $\mu_1$ is bounded by the maximal time particles spend in $\Omega$, which is uniformly bounded. This implies $\mu_1(x,t)\leq\mu_{\max}$ for some $\mu_{\max}$. Now we will prove Lipschitz continuity.

Let $t\in(0,T)$ be fixed and let $x,\tx\in\overline{\Omega}$ such that $|x-\tx|\leq\varepsilon$, where $\varepsilon$ will be chosen sufficiently small in the following. Due to the assumptions on $a$, the pathlines passing through $x$ and $\tx$ are uniquely determined and originate at some $(x_0,t_0)$ and $(\tx_0,\tt_0)$, respectively. In other words,
\begin{equation}
x=S(t;x_0,t_0),\quad \tx=S(t;\tx_0,\tt_0).
\nonumber
\end{equation}
Without loss of generality, let $t_0\ge \tt_0>0$, hence $x_0,\tx_0\in\pd\Omega^-$. The case when $\tt_0=0$ can be treated similarly. Furthermore, we assume that $x_0$ is not a vertex of $\pd\Omega^-$ -- this case will be treated at the end of the proof. 

Since $a$ is uniformly bounded and Lipschitz continuous in space, there exists $\delta$, such that if $t\in[0,T)$ and $\mathrm{dist}\{x,\pd\Omega^-\}\leq\delta$ then
\begin{equation}
-a(x,t)\cdot\bn\ge a_{\min}/2,
\label{lem:mu_reg_2}
\end{equation}
due to (\ref{lem:mu_reg_ass}). If $\varepsilon$ is sufficiently small, then for the distance of the two considered pathlines at time $t_0$ we have $|x_0-S(t_0;\tx_0,\tt_0)| =|S(t_0;x_0,t_0)-S(t_0;\tx_0,\tt_0)|\leq\delta$. Since $x_0\in\pd\Omega^-$ this means that $S(t_0;\tx_0,\tt_0)$ is in the $\delta$-neighborhood of $\pd\Omega^-$ and by (\ref{lem:mu_reg_2}), $\mathrm{dist\{S(\vartheta;\tx_0,\tt_0),\pd\Omega^-\}}$  decreases as $\vartheta$ goes from $t_0$ to $\tt_0$ with a rate of at least $a_{\min}/2$ due to the uniformity of the bound (\ref{lem:mu_reg_2}). Therefore, $S(\vartheta;\tx_0,\tt_0)$ stays in the $\delta$-neighborhood of $\pd\Omega^-$ for all $\vartheta\in[\tt_0,t_0]$ and
~\begin{equation}
-a(S(\vartheta;\tx_0,\tt_0),\vartheta)\cdot\bn\ge a_{\min}/2
\nonumber
\end{equation}
for all $\vartheta\in[\tt_0,t_0]$. Moreover, since $x_0$ lies in the interior of an edge on $\pd\Omega^-$, by choosing $\varepsilon$ sufficiently small, we can ensure that $\tx_0$ also lies on this edge (face).

Now we estimate $|\mu(x,t)-\mu(\tx,t)|=|t_0-\tt_0|$. We have
\begin{equation}
\begin{split}
x-x_0 &=\int_{t_0}^{t}\frac{\d S}{\d t}(\vartheta;x_0,t_0)\d\vartheta =\int_{t_0}^{t} a(S(\vartheta;x_0,t_0),\vartheta)\d\vartheta,\\
\tx-\tx_0 &=\int_{\tt_0}^{t}\frac{\d S}{\d t}(\vartheta;\tx_0,\tt_0)\d\vartheta =\int_{\tt_0}^{t} a(S(\vartheta;\tx_0,\tt_0),\vartheta)\d\vartheta
\nonumber
\end{split}
\end{equation}
Subtracting these two identities gives us
\begin{equation}
\begin{split}
x_0-\tx_0 =x-\tx +\int_{\tt_0}^{t_0} a(S(\vartheta;\tx_0,\tt_0),\vartheta)\d\vartheta +\int_{t_0}^{t}a(S(\vartheta;\tx_0,\tt_0),\vartheta) -a(S(\vartheta;x_0,t_0),\vartheta)\d\vartheta.
\label{lem:mu_reg_5}
\end{split}
\end{equation}
If we consider $\bn$, the normal to $\pd\Omega^-$ at $x_0$, we see that $(x_0-\tx_0)\cdot\bn=0$, as both $x_0$ and $\tx_0$ lie on the same edge on $\pd\Omega^-$. Therefore, if we multiply (\ref{lem:mu_reg_5}) by $-\bn$, we get
\begin{equation}
\begin{split}
0 =(x-\tx)\cdotp\!(-\bn) +\int_{\tt_0}^{t_0} a(S(\vartheta;\tx_0,\tt_0),\vartheta)\cdotp\!(-\bn)\d\vartheta +\int_{t_0}^{t}\big(a(S(\vartheta;\tx_0,\tt_0),\vartheta) -a(S(\vartheta;x_0,t_0),\vartheta)\big)\cdotp\!(-\bn)\d\vartheta.
\label{lem:mu_reg_6}
\end{split}
\end{equation}
Due to (\ref{lem:mu_reg_2}), we can estimate the first integral as
\begin{equation}
\int_{\tt_0}^{t_0} a(S(\vartheta;\tx_0,\tt_0),\vartheta)\cdotp\!(-\bn)\d\vartheta \ge \frac{a_{\min}}{2}(t_0-\tt_0).
\nonumber
\end{equation}
As for the second integral in (\ref{lem:mu_reg_6}), we have
\begin{equation}
\begin{split}
\Big|\int_{t_0}^{t}\big(a(S&(\vartheta;\tx_0,\tt_0),\vartheta) -a(S(\vartheta;x_0,t_0),\vartheta)\big)\cdotp\!(-\bn) \d\vartheta\Big|\\
&\leq \hat{T}L_a \sup_{\vartheta\in(t_0,t)}|S(\vartheta;\tx_0,\tt_0) -S(\vartheta;x_0,t_0)|\leq \hat{T}L_a (e^{\hat{T}L_a}-1)\|x-\tx\|,
\label{lem:mu_reg_8}
\end{split}
\end{equation}
where $L_a$ is the Lipschitz constant of $a$ with respect to $x$. The last inequality in (\ref{lem:mu_reg_8}) follows from standard results on ordinary differential equations, namely continuous dependance of the solution on the initial condition -- here we consider the ODEs defining $S(\cdotp;\cdotp,\cdotp)$ backward in time on the interval $t_0,t$ with ``initial" conditions $x$ and $\tx$ at time $t$.

Finally $|(x-\tx)\cdotp\!(-\bn)|\leq \|x-\tx\|$. Therefore, we get from (\ref{lem:mu_reg_6})
\begin{equation}
\frac{a_{\min}}{2}|\mu(x,t)-\mu(\tx,t)| =\frac{a_{\min}}{2}|t_0-\tt_0| =\frac{a_{\min}}{2}(t_0-\tt_0) \leq C|x-\tx|,
\nonumber
\end{equation}
where $C$ is independent of $x,\tx,t_0,\tt_0,t$. Dividing by $a_{\min}/2>0$ gives Lipschitz continuity of $\mu(\cdot,t)$. 

Now we return to the case when $x_0$ is a vertex of $\pd\Omega^-$. Reasoning as in the preceding, by choosing $\varepsilon$ sufficiently small, we can ensure that $\tx_0\in\pd\Omega^-$ is sufficiently close to $x_0$, i.e $\tx_0$ lies on one of the edges adjoining $x_0$. Then we can again multiply (\ref{lem:mu_reg_6}) by $-\bn$, the normal to $\pd\Omega^-$ at $\tx_0$. Hence $(\tx-\tx_0)\cdotp\!(\bn)=0$ will also be satisfied and we can proceed as in the previous case. 
\qed
\\

Having established uniform Lipschitz continuity of $\mu_1$ in space, we can prove its Lipschitz continuity with respect to time. This implies the differentiability of $\mu_1$ with respect to $x$ and $t$ a.e. in $Q_T$, therefore the left-hand side of (\ref{mu_construct:1}) is well defined a.e. in $Q_T$ and this expression is equal to the derivative of $\mu_1$ along pathlines. This is important, as all the following considerations are thus justified.

\begin{lemma}
\label{lem:mu_reg_t}
Let $a$ satisfy the assumptions of Lemma \ref{lem:mu_reg}. Then $\mu_1$ is uniformly Lipschitz continuous with respect to time $t$: there exists $L_t\ge 0$ such that
\begin{equation}
|\mu(x,t)-\mu(x,\tt)|\leq L_t|t-\tt|
\nonumber
\end{equation}
for all $x\in\Omega$ and $t,\tt\in[0,T)$.
\end{lemma}
\proof
Let e.g. $\tt>t$ and denote $(x_0,t_0)$ and $(\tx_0,\tt_0)$ denote the origin of the pathlines passing through $(x,t)$ and $(x,\tt)$, respectively. Therefore $x=S(t;x_0,t_0) =S(\tt;\tx_0,\tt_0)$. Denote $\tx=S(t;\tx_0,\tt_0)$. Then
\begin{equation}
|x-\tx| =\Big|\int_{t}^{\tt} a(S(\vartheta;\tx_0,\tt_0),\vartheta) \d\vartheta\Big| \leq \|a\|_{L^{\infty}(Q_T)}|t-\tt|.
\label{lem:mu_reg_t_1}
\end{equation}
Therefore by (\ref{mudef})
\begin{equation}
\begin{split}
|\mu(x,t)-\mu(x,\tt)| &\leq |t-\tt|+|\tt_0-t_0| =|t-\tt|+|\mu(x,t)-\mu(\tx,t)| \leq |t-\tt|+L_\mu|x-\tx|\\
&\leq (1+L_\mu \|a\|_{L^{\infty}(Q_T)})|t-\tt|
\nonumber
\end{split}
\end{equation}
due to (\ref{lem:mu_reg_t_1}) and Lemma \ref{lem:mu_reg}. This completes the proof.

\qed

\noindent\textbf{Example:} In simple cases, the function $\mu_1$ can be explicitly written down. As a trivial example, we take the one-dimensional stationary flow field $a(x,t)=x+1$ on $\Omega=(0,1)$ and the time interval $(0,+\infty)$. Then (\ref{mu_construct:2}) can be easily solved to obtain
\begin{equation}
S(t;x_0,t_0)=(x_0+1)e^{t-t_0}-1.
\nonumber
\end{equation}
The $(x,t)$--plane is then separated into two regions separated by the pathline $S(t;0,0)$, which is the curve $x=e^t-1$, i.e. $t=\ln(x+1)$. For points $(x,t)$ beneath this curve, i.e. $t\leq\ln(x+1)$, we have $t_0=0$, hence $\mu_1(x,t)=t-t_0=t$. For points $(x,t)$ above the separation curve, we have $x_0=0$, hence
\begin{equation}
x=S(t;0,t_0)=e^{t-t_0}-1 \quad\Longrightarrow \quad \mu_1(x,t)=t-t_0=\ln(x+1).
\nonumber
\end{equation}
Altogether, we have
\begin{equation}
\mu_1(x,t)=\begin{cases}   
t & \text{if $t\leq\ln(x+1)$,} \\
\ln(x+1)& \text{otherwise.}
\end{cases}
\nonumber
\end{equation}
This function is globally bounded and globally Lipschitz continuous. We note that the standard exponential scaling trick which gives exponential growth of the error corresponds to taking $\mu_1(x,t)=t$ for all $x,t$, which is an unbounded function.

\section{Conclusion and future work}
In this paper we have derived a priori error estimates for a linear advection reaction problem with inlet and outlet boundary conditions in the $L^\infty(L^2)$ and $L^2(L^2)$ norms. Unlike previous works, the analysis was performed without the usual ellipticity assumption $c-\tfrac{1}{2}\diver a\ge 0$. This is achieved by applying a general exponential scaling transformation in space and time to the exact and discrete solutions of the problem. We considered the case when the time spent by particles carried by the flow field inside the spatial domain $\Omega$ is uniformly bounded by some $\widehat{T}$. The resulting error estimates are of the order $Ch^{p+1/2}$, where $C$ depends exponentially on $\widehat{T}$ (which is a constant) and not on the final time $T$, as would be expected from the use of Gronwall's inequality. Effectively, due to the exponential scaling, we apply Gronwall's lemma in the Lagrangian setting along pathlines, which exist only for time at most $\widehat{T}$, and not in the usual Eulerian sense.

As for future work, we plan to extend the analysis also to fully discrete DG schemes with discretization in time. Furthermore, we wish to extend the analysis to nonlinear convective problems, following the arguments of \cite{ZhangShu} and \cite{kucera} to obtain error estimates without the exponential dependence on time of the error estimates.

\section*{Acknowledgements}
The work of V. Ku\v{c}era was supported by the J. William Fulbright Commission in the Czech Republic and research project No. 17-01747S of the Czech Science Foundation. The work of C.-W. Shu was supported by DOE grant DE-FG02-08ER25863
and NSF grant DMS-1418750.



\end{document}